\RequirePackage{fix-cm}
\documentclass[smallextended]{svjour3}       
\smartqed  
\usepackage{graphicx}
\usepackage{amssymb}
\usepackage{amsmath}
\usepackage{amscd}
\newcommand{\be}{\begin{equation}}
\newcommand{\ee}{\end{equation}}
\newcommand{\ba}{\begin{aligned}}
\newcommand{\ea}{\end{aligned}}
\begin{document}

\title{A coupling approach to Doob's theorem
}


\author{Alexei   Kulik     \and
        Michael Scheutzow 
}


\institute{A. Kulik \at
              Institute of Mathematics, NAS of Ukraine, 3, Tereshchenkivska str., 01601  Kyiv, Ukraine \\
              Tel.: +38-04-42793994\\
              \email{kulik.alex.m@gmail.com}     
           \and
           M. Scheutzow \at
              Institut f\"ur Mathematik, MA 7-5, Technische Universit\"at Berlin, Str. des 17. Juni 136, 10623 Berlin, Germany\\
              Tel.: +49-30-31425767\\
              Fax: +49-30-31421695\\
              \email{ms@math.tu-berlin.de}
}

\date{Received: date / Accepted: date}

\maketitle

\begin{abstract}
We provide a coupling proof of Doob's theorem which says that the transition probabilities of a regular Markov process which has an invariant probability measure $\mu$
converge to $\mu$ in the total variation distance. In addition we show that non-singularity (rather than equivalence) of the transition probabilities suffices to
ensure convergence of the transition probabilities for $\mu$-almost all initial conditions.

\end{abstract}

\section{Introduction}
\label{intro}
Doob's theorem, as formulated in \cite{DZ96}, p.43, states that if the stochastically continuous Markov semigroup $P_t,\,t\ge 0$ with Polish state space $(E,d)$ has an
invariant probability measure (ipm)  $\mu$ and is $t_0$-{\em regular} for some $t_0>0$, then $\mu$ is unique and all transition probabilities converge to $\mu$ in the
total variation distance. Here, $t_0$-regular means that all transition probabilities $P_{t_0}(x,.)$ are mutually equivalent. One common way to check  $t_0$-regularity is
to show that the Markov semigroup is irreducible and strong Feller, see \cite{DZ96}, Proposition 4.1.1 (known as Khas'minskii's theorem).

In fact, Da Prato and Zabczyk formulate and prove Doob's theorem with respect to strong convergence (which is weaker than total variation convergence) and refer
the reader to \cite{St94} and \cite{Se97} for different proofs of total variation convergence. Neither of the proofs is short and elementary. In particular,
neither of the proofs uses {\em coupling} which has been a powerful tool to prove convergence of transition probabilities in the past decades. The aim of this article
is to provide a coupling proof of Doob's theorem. At the same time we generalize the result in various directions: instead of a Polish state space, we just require a mild
condition on the measurable space $(E,{\mathcal{E}})$. Further, we allow the infimum of the times $t>0$ for which $P_t(x,.)$ and $P_t(y,.)$ are equivalent to depend on the pair $(x,y)$
without being uniformly bounded from above. We also show that we can replace equivalance of the transition probabilities by the much weaker property of non-singularity but in this case
convergence of the transition probabilities only holds for $\mu$-almost all initial conditions $x$. In the next section, we formulate the main results in the discrete time setting
(in Remark  \ref{Bemerkung} we say why the continuous time claim follows)  and provide an example showing that the weaker assumption in Corollary \ref{cor1} does not guarantee the
conclusion of Theorem \ref{thm1}.

\section{Main results}
\label{sec:1}

Let $X_n, n\in {\mathbb{Z}}_+$ be a Markov chain with the state space $(E, \mathcal{E})$. The measurable space $(E, \mathcal{E})$ is assumed to be
countably generated. We also assume that the diagonal $\Delta=\{(x,x), x\in E\}$ belongs to $\mathcal{E}\otimes \mathcal{E}$.
A typical example of such a space is a Borel measurable space, e.g. a Polish space $E$ endowed with the Borel $\sigma$-algebra $\mathcal{E}$.

Transition probabilities  and $n$-step transition probabilities for $X$ are denoted respectively by $P(x,dy)$ and $P_n(x,dy)$.
The law of the sequence $\{X_n\}$ in $(E^\infty, \mathcal{E}^{\otimes \infty})$ with initial distribution $\mathrm{Law}\, (X_0)=\mu$
is denoted by $\mathbb{P}_\mu$, the respective expectation is denoted by $\mathbb{E}_\mu$; in case $\mu=\delta_x$ we write
simply $\mathbb{P}_x, \mathbb{E}_x$.

Recall that an invariant probability measure for $X$ is a probability measure $\mu$ on $(E, \mathcal{E})$ such that
\begin{equation}\label{invar}
\mu(dy)=\int_EP(x,dy)\mu(dx).
\end{equation}
Equivalently, a probability measure $\mu$ is invariant if the sequence $\{X_n, n\in \mathbb{Z}_+\}$ is strictly stationary under $\mathbb{P}_\mu$.

We use the usual relations for probability measures $\mu, \nu$ on $(E, \mathcal{E})$: $\mu$ and $\nu$ are equivalent (notation $\mu\sim\nu$)
if each of them is absolutely continuous w.r.t. the other; $\mu$ and $\nu$ are singular (notation $\mu\perp\nu$) if there exists $A\in \mathcal{E}$
such that $\mu(A)=1, \nu(A)=0$; otherwise  $\mu$ and $\nu$ are non-singular (notation $\mu\not\perp\nu$). The total variation distance between
probability measures $\mu, \nu$ on $(E, \mathcal{E})$ is the total variation of the signed measure $\mu-\nu$ (notation $\|\mu-\nu\|$).

\begin{theorem}\label{thm1} Assume that for each $x,y\in E$  there exists $n=n_{x,y}$ such that
\be\label{ass1}
P_n(x, \cdot)\sim P_n(y, \cdot).
\ee

Then there exists at most one ipm for the chain $X$. If an ipm $\mu$ exists, then for every $x\in E$
\be\label{conv}
\|P_n(x, \cdot)-\mu\|\to 0, \quad n\to \infty.
\ee
\end{theorem}
The proof of Theorem \ref{thm1} is given in Section \ref{sec:2} below. The following theorem shows to what extent the basic assumption (\ref{ass1})  can be relaxed.

\begin{theorem}\label{thm2} Let $X$ be a Markov chain which has an ipm $\mu$. Assume further that for $\mu\otimes\mu$-almost all $(x,y)\in E\times E$  there exists $n=n_{x,y}$ such that
\be\label{ass2}
P_n(x, \cdot)\not\perp P_n(y, \cdot).
\ee

Then (\ref{conv}) holds true for $\mu$-almost all $x\in E$.
\end{theorem}

The following corollary is simple and straightforward.
\begin{corollary}\label{cor1} Assume that for each $x,y\in E$  there exists $n=n_{x,y}$ such that (\ref{ass2}) holds true.

Then there exists at most one ipm for the chain $X$. If an ipm $\mu$ exists, then (\ref{conv}) holds true for $\mu$-a.a. $x\in E$.
\end{corollary}

To prove uniqueness of an ipm (which is  the only addition to Theorem \ref{thm2}), let us assume that there exist two different ipm's
$\mu_1, \mu_2$, and consider the averaged ipm $\mu=(1/2)(\mu_1+\mu_2)$. Applying Theorem \ref{thm2} first to $\mu_1$ and then to $\mu$, we get a contradiction: because $\mu_1$ is absolutely continuous w.r.t. $\mu$, we get that for $\mu_1$-a.a. $x\in E$ the transition probabilities $P_n(x, \cdot)$ converge both to $\mu_1$ and to $\mu$, but $\mu_1\not=\mu$. \qed

\begin{remark} Note that the condition of Theorem \ref{thm2} alone does not yield uniqueness of the ipm for $X$: a simple  counter-example is given by a chain with a finite state
space with at least two mutually disconnected classes of states.
\end{remark}

\begin{remark}\label{Bemerkung}
Note that the results of Theorem \ref{thm1}, Theorem \ref{thm2}, and Corollary \ref{cor1} are also true in the continuous time case when $P_t,\, t\ge 0$ is a
Markov semigroup (no regularity in $t$ is required). To see this, note that uniqueness of an ipm $\mu$ for the discretized chain
$P_n,\, n \in {\mathbb{N}}_0$ implies uniqueness of an ipm for $P_t,\,t \ge 0$ and
that for any (discrete of continuous time) Markov semigroup $(P_r)$ with ipm $\mu$, the function $r \mapsto \|P_r(x,.)-\mu\|$
is non-increasing.
\end{remark}

\begin{example}
The following example shows that the assumptions of Corollary  \ref{cor1} do not imply the conclusion of Theorem \ref{thm1}. Equip $E={\mathbb{Z}_+}$ with the
discrete $\sigma$-algebra $\mathcal{E}$ and define the transition probabilities by $p_{0,0}=1$, $p_{i,i-1}=1/3$ and $p_{i,i+1}=2/3$ for $i \ge 1$. Then $\mu=\delta_0$ is the unique invariant
probability measure, the assumptions of Corollary \ref{cor1} hold but $P_n(i,.)$ does not converge to $\mu$ for any $i \neq 0$.
\end{example}

\begin{remark}
Note that in the discrete case (i.e.~$E$ is finite or countably infinite) the assumptions in  Theorem \ref{thm2} and in Corollary \ref{cor1}  are also
necessary for the respective conclusion to hold. This is not true for Theorem \ref{thm1} however
as the example $E={\mathbb{Z_+}}$ with $p_{0,0}=p_{0,1}=1/2$, $p_{i,i-1}=2/3$ and $p_{i,i+1}=1/3$ for $i \ge 1$ shows.
\end{remark}

\section{Proofs of Theorem \ref{thm1} and Theorem \ref{thm2}}\label{sec:2}
 \paragraph{An auxiliary construction.}
 Denote for $N\in \mathbb{N}, p\in (0,1)$
 $$
 C_{N,p}=\{(x,y): \|P_N(x, \cdot)-P_N(y, \cdot)\|\leq 2(1-p)\}.
 $$
 Since
 $$
 \nu\not \perp\bar \nu \Leftrightarrow \|\nu-\bar\nu\|<2,
 $$
 the assumption of  Theorem \ref{thm2} (and therefore the stronger assumption of  Theorem \ref{thm1}, as well) yield that there exist
$N\in \mathbb{N}, p\in (0,1)$ such that
$$
 (\mu\otimes \mu)(C_{N,p})>0.
$$
 In the sequel we fix these values $N,p$ and write simply $C$ instead of $C_{N,p}$. In addition we assume that $N=1$. Remark \ref{Bemerkung} shows
that this is no loss of generality.
%

  \paragraph{An outline of the method.} Our aim, in fact, is to prove the convergence
\be\label{conv_point}
\|P_n(x_1, \cdot)-P_n(x_2, \cdot)\|\to 0, \quad n\to \infty
\ee
either for all $(x_1,x_2)\in E\times E$ in the case considered in Theorem \ref{thm1}, or for $\mu\otimes \mu$-a.a. $(x_1,x_2)\in E\times E$ in the case considered in Theorem \ref{thm2}. Once (\ref{conv_point}) is proved, the required convergence (\ref{conv}) follows
using the representation \eqref{invar} and the triangle inequality.
The following fact  is  well-known (\cite{Thor}, p.14): for  any  two random elements $\xi_1, \xi_2$, defined on a same probability space $(\Omega, \mathcal{F}, \mathbb{P})$,   valued in $(E, \mathcal{E})$, and such that  $\mathrm{Law}(\xi_i)=\nu_i, i=1,2$,  one has
\be\label{coup}
\|\nu_1-\nu_2\|\leq 2{\mathbb{P}}(\xi_1\not=\xi_2).
\ee
Hence for any sequence  $\{Z_n=(Z_n^1, Z_n^2), n\in \mathbb{Z}_+\}$ such that the laws of
$\{Z_n^i, n\in \mathbb{Z}_+\}, i=1,2$ equal respectively $\mathbb{P}_{x_i}, i=1,2$, one has a bound
$$
\|P_n(x_1, \cdot)-P_n(x_2, \cdot)\|\leq 2{\mathbb{P}}(Z_n^1\not=Z_n^2).
$$
This is the essence of the famous \emph{coupling} approach which dates back to W.~D\"oblin \cite{Do40}: to prove (\ref{conv_point}) one should construct a sequence $Z$ which verifies the above assumption (such a sequence is usually called a \emph{coupling} for $X$) in such a way that
\be\label{conv_prob}
{\mathbb{P}}(Z_n^1\not=Z_n^2)\to 0, \quad n\to \infty.
\ee
\paragraph{Construction of the coupling.} The sequence $Z$ will be taken as a Markov chain on $E\times E$ with suitably constructed transition probability. The first part of this construction is based on the fact that a proper choice of the pair $(\xi_1,\xi_2)$ may turn inequality
(\ref{coup}) into an identity. This fact, sometimes called the Coupling Lemma, is well known; the law of any such pair $(\xi_1,\xi_2)$ is called a \emph{maximal coupling}. We refer to \cite{Thor} Section 1.4,  where the construction of a maximal coupling based on the \emph{splitting representation} of a random variable is given. In our framework we use essentially the same construction, but with modifications which are caused by the necessity (a) to deal with transition probabilities instead of measures, and therefore to take care of measurability issues; (b) to manage properly the law of the pair ``outside of the diagonal''. Namely, we have the following statement (the proof is given in the Appendix).

\begin{lemma}\label{lem1} (The Coupling Lemma for transition probabilities) Let $(E, \mathcal{E})$ be a countably generated measurable space.

Then for any Markov kernel $P(x, dy)$ on $(E, \mathcal{E})$ there exists a Markov kernel $Q((x_1,x_2), dy_1dy_2)$ on $(E\times E, \mathcal{E}\otimes \mathcal{E})$ such that for every $(x_1,x_2)\in E\times E$:
\begin{itemize}
  \item[(i)] $Q((x_1,x_2), \Delta)=1-(1/2)\|P(x_1, \cdot)-P(x_2, \cdot)\|$;
  \item[(ii)] the measure $Q((x_1,x_2), dy_1dy_2)$ restricted to $(E\times E)\setminus \Delta$ is absolutely continuous w.r.t. $P(x_1, dy_1)\otimes P(x_2, dy_2)$.
\end{itemize}
\end{lemma}

\begin{remark} We refer a reader to \cite{Ver} for another version of the Coupling Lemma which takes into account measurability issues; the  measurable space $(E, \mathcal{E})$ therein is assumed to be a Borel one.
\end{remark}

Denote
$$
R((x_1,x_2), dy_1dy_2)=P(x_1, dy_1)\otimes P(x_2, dy_2),
$$
which is just the transition probability of the Markov chain in $E\times E$ whose components are independent and each of the components is a Markov chain with the transition probability $P(x,dy)$. Such a chain is usually called an \emph{independent coupling}, and below we denote it by $W=\{W_n, n\in \mathbb{Z}_+\}.$

Finally, we define the transition probability for $Z$ by
$$
S\big((x_1,x_2),.\big):=\left\{
\begin{array}{ll}
Q\big( (x_1,x_2),.\big)&\mbox{if } (x_1,x_2) \in C\\
R\big( (x_1,x_2),.\big)&\mbox{otherwise,}
\end{array}
\right.
$$
with the set $C$ defined at the beginning of the proof. By construction, $Z$ is a coupling for $X$, and our aim is to prove (\ref{conv_prob}) with $Z_0=(x_1, x_2)$ either for all $(x_1, x_2)$ in the case of Theorem \ref{thm1}, or for $\mu\otimes\mu$-a.a. $(x_1, x_2)$ in the case of Theorem \ref{thm2}.

\paragraph{Proof of (\ref{conv_prob}).} Note that $\Delta\subset C$, and for any point $(x,x)\in \Delta$
$$
Q((x,x), \Delta)=1.
$$
Hence, by the construction the following property holds: once $Z$ hits $\Delta$, all the subsequent values of $Z$ a.s. stay in $\Delta$ (``once the components are coupled they stay coupled''). Therefore
$$
{\mathbb{P}}(Z_n^1\not=Z_n^2)={\mathbb{P}}(T>n), \quad T=\inf\{m: Z_m^1\not=Z_m^2\},
$$
with the usual convention $\inf\emptyset=\infty$.

Consider the sequence of stopping times
$$
\tau_0=0, \quad \tau_k=\inf\{n>\tau_{k-1}: Z_n\in C\}, \quad k\geq 1,
$$
and \emph{assume} for a moment that we know that
\be\label{rec_tau}
\tau_k<\infty, \quad k\geq 1
\ee
with probability 1. Clearly, for any $k\geq 1$
$$
\{T>\tau_k\}=\bigcup_{n=1}^\infty\{\tau_k=n, Z_1\not \in \Delta, \dots, Z_n\not\in \Delta\}\in \mathcal{F}_{\tau_k};
$$
here  $\{\mathcal{F}_n\}$ denotes the natural filtration of the sequence $Z$. By the construction of $Z$,
$$
{\mathbb{P}}(Z_{\tau_k+1}\in \Delta|\mathcal{F}_{\tau_k})\quad \begin{cases}\geq p, & T>\tau_k,\\
=1, &\hbox{otherwise}. \end{cases}
$$
Hence for the sequence $p_k={\mathbb{P}}(T\leq \tau_k), k\geq 1$ one has
$$
p_{k+1}\geq p(1-p_k)+p_k, \quad k\geq 1,
$$
and therefore $p_k\to 1, k\to \infty$, which together with (\ref{rec_tau}) yields (\ref{conv_prob}).

\paragraph{Proof of (\ref{rec_tau}): the Recurrence Lemma.} We have reached the last and the crucial step in the proof: we need to prove that the set $C$, which in a sense is ``favorable for the subsequent coupling attempt'', is a.s. visited by $Z$ infinitely often. We use the following lemma, whose proof is given in the Appendix.

\begin{lemma}\label{lem2} (The Recurrence Lemma) Assume that the  Markov chain $X$ satisfies the  condition of Theorem \ref{thm2}.

Then for any $B\in \mathcal{E}$ with $\mu(B)>0$, for $\mu$-a.a. $x\in E$
\be\label{rec}
\mathbb{P}_x(X_n\in B\hbox{ infinitely often})=1.
\ee

If, in addition, the condition of Theorem \ref{thm1} holds true, then (\ref{rec}) holds true for every $x\in E$.
\end{lemma}

Now we can finish the whole proof; consider first the case of Theorem \ref{thm1}.  The independent coupling $W$ verifies the assumptions of Lemma \ref{lem2} with $E\times E$ instead of $E$ and $\mu\otimes\mu$ instead of $\mu$. Because $(\mu\otimes\mu)(C)>0$, this yields
\be\label{31}
{\mathbb{P}}(W_n\in C\hbox{ infinitely often})=1
\ee
 for all initial values $W_0=(x_1, x_2)\in E\times E$.

 Observe that  up to $\tau_1$ the law of $Z$ coincides with the law of the independent coupling $W$ up to its first visit to $C$, hence by (\ref{31})
 $$
{\mathbb{P}}(\tau_1<\infty)=1
 $$
for all initial values $Z_0=(x_1, x_2)\in E\times E$.  Hence $Z$ a.s. performs at least one ``coupling attempt''. If this attempt is successful, i.e. $T\leq \tau_1+1$, then a.s.
$\tau_2=\tau_1+1, \tau_3=\tau_2+1, \dots$ because $\Delta\subset C$ and, once the components of $Z$ are coupled, they stay coupled. In that case (\ref{rec_tau}) holds true. If ``the first coupling attempt is not successful'', the chain $Z$ afterwards again performs as the independent coupling $W$ up to the time moment $\tau_2$. Applying Lemma \ref{lem2} once again and the strong Markov property of $Z$, we get
$$
{\mathbb{P}}(\tau_2<\infty)=1.
$$
Iterating this argument, we get  (\ref{rec_tau}) for all initial values $Z_0=(x_1, x_2)\in E\times E$, which completes the proof of Theorem \ref{thm1}.

Let us proceed with Theorem \ref{thm2}; in that case it is convenient to prove (\ref{rec_tau}) for a  version of $Z$ with $\mathrm{Law}(Z_0)=\mu\otimes \mu$.    The argument, completely the same as above, proves that $\tau_1<\infty$ a.s., and if ``the first coupling attempt is  successful'' then (\ref{rec_tau}) holds true. Consider the law of $Z_{\tau_{1}+1}$ conditioned by the event that ``the first coupling attempt is not successful''; that is, $\{Z_{\tau_{1}+1}\not \in \Delta\}$. By the choice of the law of $Z_0$ and the construction of the kernel $Q$, this law is absolutely continuous w.r.t. $\mu\otimes\mu$. %
 Using this and  the strong Markov property of $Z$ at the stopping time $\tau_1+1$, we apply Lemma \ref{lem2} once again and get
$$
 {\mathbb{P}}(\tau_2<\infty)=1.
$$
Iterating this argument, we get  (\ref{rec_tau}) for $Z$ with  $\mathrm{Law}(Z_0)=\mu\otimes \mu$, which completes the proof of Theorem \ref{thm2}.\qed

\appendix
\section{Proofs of the auxiliary lemmas}
\subsection{Proof of Lemma \ref{lem1}} Denote
$$
\Lambda(x_1,x_2; dy)=(1/2)\Big(P(x_1, dy)+P(x_2,dy)\Big).
$$
Then for every $x_1, x_2$ $P_i(x_i, dy)\ll \Lambda(x_1,x_2; dy), i=1,2$. Let us show that respective Radon-Nikodym derivatives can be chosen in a jointly measurable way; that is, there exist measurable functions $f_i:E\times E\times E\to \mathbb{R}^+, i=1,2$ such that
$$
P(x_i, A)=\int_{A}f_i(x_1, x_2,y) \Lambda(x_1,x_2; dy), \quad i=1,2, \quad x_1, x_2\in E, \quad A\in \mathcal{E}.
$$
Let $\mathcal{E}_0$ be a countable algebra which generates $\mathcal{E}$, then the countable class $\mathcal{H}$, which consists of all functions representable in the form of a finite sum $\sum_kc_k1_{A_k} , \{c_k\}\subset \mathbb{Q}, \{A_k\}\subset \mathcal{E}_0$, is dense in $L_1(E,\lambda)$ for any probability measure $\lambda$ on $(E, \mathcal{E})$.

Consider a Radon-Nikodym derivative
$$
\rho_{x_1,x_2}^1(y)={P(x_1, dy)\over \Lambda(x_1,x_2; dy)},
$$
then for every $\varepsilon>0$ there exists $h\in \mathcal{H}$ such that
$$
\sup_{A\in \mathcal{E}_0}\int_A(\rho_{x_1,x_2}^1(y)-h(y))\Lambda(x_1,x_2; dy)<{\varepsilon\over 2}.
$$
Observe that this relation is equivalent to
\be\label{A1}
\sup_{A\in \mathcal{E}_0}\left(P(x_1, A)-\int_A h(y)\Lambda(x_1,x_2; dy)\right)<{\varepsilon\over 2},
\ee
and yields
$$
\int_E|\rho_{x_1,x_2}(y)-h(y)|\Lambda(x_1,x_2; dy)<\varepsilon.
$$
Fix some enumeration of the class $\mathcal{H}=\{h_m, m\in \mathbb{N}\}$, and denote for $n\geq 1$  by $m(x_1, x_2; n)$ the minimal $m\geq 1$ such that (\ref{A1}) holds true for $h=h_m$ with $\varepsilon=2^{-n-1}$. Then $m(\cdot;n):E\times E\to \mathbb{N}$ is measurable, and therefore
$$
f_1^n(x_1,x_2,y)=h_{m(x_1,x_2;n)}(y)
$$
is measurable as a function $E\times E\times E\to \mathbb{R}$. When $x_1,x_2$ are fixed, the sequence $\{f_1^n(x_1,x_2,y), n\geq 1\}$ converges to $\rho_{x_1,x_2}^1(y)$ for $\Lambda(x_1,x_2;\cdot)$-a.a. $y$: this follows from the Borel-Cantelli lemma because, by the construction,
$$
\|f_1^n(x_1,x_2,\cdot)-\rho_{x_1,x_2}^1\|_{L_1(E, \Lambda(x_1,x_2;\cdot))}\leq 2^{-n}.
$$
Therefore the function
$$
f_1(x_1,x_2,y)=\begin{cases}\lim_{n\to \infty}f_1^n(x_1,x_2,y),&\hbox{if the limit exists}\\
0,&\hbox{otherwise}\end{cases}
$$
gives the required measurable version of the Radon-Nikodym derivative for $P(x_1, dy)$ (the construction for $P(x_2, dy)$ is the same).

We finish the proof by repeating essentially the construction  from \cite{Thor} Section 1.4, based on the splitting representation for probability laws. Write
$$
g(x_1, x_2,y)=\min\Big(f_1(x_1,x_2,y), f_2(x_1,x_2,y)\Big), \quad p(x_1, x_2)=\int_{E}g(x_1, x_2,y)\Lambda(x_1,x_2;dy),
$$
$$
\Theta(x_1,x_2;dy)={g(x_1, x_2,y)\over p(x_1, x_2)}\Lambda(x_1, x_2;dy)
$$
with the convention that (anything)$/0=1$. Then we have representations
$$
P(x_i,dy)=p(x_1,x_2)\Theta(x_1,x_2;dy)+(1-p(x_1,x_2))\Sigma_i(x_1,x_2;dy), \quad i=1,2
$$
with probability kernels $\Theta, \Sigma_1, \Sigma_2$ and measurable $p:E\times E\to [0,1]$.

Note that the mapping $E\ni x\mapsto (x,x)\in E\times E$ is  $\mathcal{E}- \mathcal{E}\otimes \mathcal{E}$ measurable, and denote by
$Q_1((x_1,x_2), dy_1dy_2)$ the image of $\Theta(x_1,x_2;dy)$ under this mapping. Denote also by $Q_2((x_1,x_2), dy_1dy_2)$ the product of the measures $\Sigma_i(x_1,x_2;dy_i), i=1,2$. Then
$$
Q((x_1,x_2), \cdot)=p(x_1,x_2)Q_1((x_1,x_2), \cdot)+(1-p(x_1,x_2))Q_2((x_1,x_2), \cdot)
$$
is the required kernel; observe that the assertion (ii) now holds true because $Q_2((x_1,x_2), \cdot)$ is chosen to be a product measure with the components $\Sigma_i(x_1,x_2;\cdot)\ll P(x_i, \cdot), i=1,2$.\qed

\subsection{Proof of Lemma \ref{lem2}} Denote
$$
\psi(x)=\mathbb{P}_x(X_n\in B\hbox{ infinitely often}),
$$
and consider a stationary version of $X$ with $\mathrm{Law}(X_0)=\mu$. Then the sequence $\{\psi(X_n)\}$ is stationary. But, in addition, this sequence is a L\'evy martingale:
by the Markov property of $X$, we have with probability 1
 $$
 \psi(X_n)={\mathbb{E}}_{X_n}1_{X_k\in B\mbox{ i.o.}}={\mathbb{E}}[1_{X_k\in B\mbox{ i.o.}, k\geq n}|\mathcal{F}_n]={\mathbb{E}}[1_{X_k\in B\mbox{ i.o.}, k\geq 0}|\mathcal{F}_n].
 $$
Then with probability 1
$$
 \psi(X_n)\to 1_{X_k\in B\mbox{ i.o.}, k\geq 0}, \quad n\to \infty,
$$
and hence by stationarity of $\{\psi(X_n)\}$ we have  $\psi(x)=0$ or 1 for $\mu$-a.a. $x\in E$. Denote
$$
\Psi_0=\{x:\psi(x)=0\}, \quad \Psi_1=\{x:\psi(x)=1\},
$$
and observe that because $\psi(x)\in [0,1]$ in any case, by the martingale property of $\{\psi(X_n)\}$ one has for any $n\geq 1$
$$
P_n(x,\Psi_0)=1\hbox{ for $\mu$-a.a. $x\in \Psi_0$ and } P_n(x,\Psi_1)=1\hbox{ for $\mu$-a.a. $x\in \Psi_1$.}
$$
Because of the assumption of the lemma (or Theorem \ref{thm2}), it is impossible that both $\Psi_0$ and $\Psi_1$ have positive measure $\mu$. The identity $\mu(\Psi_1)=0$
contradicts Birkhoff's ergodic theorem: with probability 1 we have
$$
{1\over N}\sum_{n=1}^N 1_B(X_n) \to \eta,\quad N\to \infty
$$ where $E \eta=\mu(B)>0$. Therefore $\mu(\Psi_0)=0,$ which implies the required identity $\mu(\Psi_1)=1$.

If, in addition, the condition of Theorem \ref{thm1} holds true, then it follows from what we have just  proved  that for any $x\in E$ there exists $n=n_x$ such that $P_n(x,\Psi_1)=1$. By the Markov property of $X$ and the definition of $\psi$, this implies $\psi(x)=1$. \qed




\end{document}